\documentclass[a4paper,reqno]{amsart}

\title{Dynamics of Mandelbrot cascades}

\author[Barral, Peyri\`ere and WEN]
{Julien Barral$^{*}$ \and Jacques Peyri\`ere$^{\dag\ddag}$ \and
Zhi-Ying WEN$^{\dag\&}$}
\thanks{$^*$ INRIA Rocquencourt, B.P. 105, 78153 Le Chesnay Cedex,
France. \texttt{Julien.Barral@inria.fr}}
\thanks{$^\dag$ Dept.\ Math., Tsinghua
University, Beijing, 100084, P.\ R.\ China. \texttt{wenzy@mail.tsinghua.edu.cn}}
\thanks{$^\ddag$ Universit\'e Paris-Sud, Math\'ematique b\^at.\ 425,
CNRS UMR 8628, 91405 Orsay Cedex, France. \texttt{Jacques.Peyriere@math.u-psud.fr}}
\thanks{\& Supported by the National Basic Research Program of China
(973 Program) (2007CB814800).}

\usepackage{amsmath,amsfonts,amssymb,mathrsfs}
\newtheorem{theorem}{Theorem}
\newtheorem{lemma}[theorem]{Lemma}
\newtheorem{proposition}[theorem]{Proposition}
\newtheorem{corollary}[theorem]{Corollary}
\theoremstyle{definition}
\theoremstyle{remark}
\newtheorem{remark}{Remark}
\newcommand{\dif}{\mathrm{d}}
\newcommand{\moment}{\mathbf{m}}
\newcommand{\T}{\mathsf{T}}
\newcommand{\s}{\mathsf{S}}
\newcommand{\Id}{\mathrm{Id}}
\newcommand{\ind}{{\mathbf{1}}}
\newcommand{\alphabet}{\mathscr{A}}
\newcommand{\mat}[1]{\mathsf{#1}} 
\DeclareMathOperator{\esp}{\mathbb{E}} 
\begin{document}

\begin{abstract}
Mandelbrot multiplicative cascades provide a construction of a
dynamical system on a set of probability measures defined by
inequalities on moments. To be more specific, beyond the first
iteration, the trajectories take values in the set of fixed points of
smoothing transformations (i.e., some generalized stable laws).

Studying this system leads to a central limit theorem and to its
functional version. The limit Gaussian process can also be obtained as
limit of an `additive cascade' of independent normal variables.
\end{abstract}

\maketitle

\section{A dynamical system}\label{dynam}

Consider the set~$\alphabet=\{0,\dots,b-1\}$, where $b\ge 2$. Set
$\alphabet^*=\bigcup_{n\ge 0}\alphabet^n$, where, by convention,
$\alphabet^0$ is the singleton~$\{\epsilon\}$ whose the only element
is the empty word~$\epsilon$. If~$w\in\alphabet^*$, we denote by~$|w|$
the integer such that $w\in \alphabet^{|w|}$. If $n\ge 1$ and
$w=w_1\cdots w_n\in \alphabet^n$ then for $1\le k\le n$ the word
$w_1\cdots w_k$ is denoted $w|_k$. By convention, $w|_0=\epsilon$.

Given~$v$ and~$w$ in~$\alphabet^n $, $v\wedge w$ is defined to be the
longest prefix common to both~$v$ and~$w$, i.e., $v|_{n_{0}}$, where
$n_{0}=\sup \{0\le k\le n\ :\ x|_n=y|_n\}$.

Let $\alphabet^\omega$ stand for the set of infinite sequences
$w=w_1w_2\cdots$ of elements of~$\alphabet$. Also, for $x\in
\alphabet^\omega $ and $n\geq 0$, let $x|_n$ be the projection of $x$
on $ \alphabet^{n}$.

If $w\in \alphabet^*$, we consider the cylinder~$[w]$ consisting of
infinite words in~$\alphabet^\omega$ whose $w$ is a prefix.

We index the closed $b$-adic subintervals of $[0,1]$ by $\alphabet^*$:
for~$w\in \alphabet^*$, we set
$$I_w =\left [\sum_{1\le k\le |w|}w_kb^{-k},
\sum_{1\le k\le |w|}w_kb^{-k}+b^{-|w|}\right ].
$$
\medskip

If $f\ :\ [0,1]\mapsto\mathbb{R}$ is bounded, for every sub-interval
$I=[\alpha,\beta]$ of $[0,1]$, we denote by~$\Delta(f,I)$ the
increment $f(\beta)-f(\alpha)$ of~$f$ over the interval~$I$.
\medskip

Let $\mathcal{P}$ the set of Borel probability measures on
$\mathbb{R}_+$. If $\mu\in\mathcal{P}$ and $p> 0$, we denote by
$\moment_p(\mu)$ the moment of order $p$ of $\mu$, i.e.,
$$
\moment_p(\mu)=\int_{\mathbb{R}_+}x^p \,\mu(\dif x).
$$ Then let $\mathcal{P}_1$ be the subset of $\mathcal{P}$ whose
elements have their first moment equal to~1:
$$
\mathcal{P}_1=\{\mu\in\mathcal{P}\ :\ \moment_1(\mu)=1\}.
$$

The smoothing transformation~$\s_\mu$ associated with~$\mu\in
{\mathcal P}$ is the mapping from~$\mathcal{P}$ to itself so defined:
If $\nu\in{\mathcal P}$, one considers~$2b$ independent random
variables, $Y(0),\, Y(1),\dots,\,Y(b-1)$, whose common probability
distribution is~$\nu$, and $W(0)$,\\$W(1),\dots,$ $\,W(b-1)$ whose
common probability distribution is~$\mu$; then $\s_\mu\nu$ is the
probability distribution of $\displaystyle b^{-1} \sum_{0\le j<b}
W(j)\,Y(j)$.

This transformation and its fixed points have been considered in
several contexts, in particular by B.~Mandelbrot who introduced it to
construct a model for turbulence and intermittence (see
\cite{M1,M2,K1,P1,KP,DL,Gu}).

In this latter case, the measure $\mu$ is in $\mathcal{P}_1$ so that
$\s_\mu$ maps $\mathcal{P}_1$ into itself.  It is known that the
condition $\int x\log(x) \,\mu(\dif x)<\log b$ is necessary and
sufficient for the weak convergence of the sequence~$\s_\mu^n\delta_1$
(where~$\delta_1$ stands for the Dirac mass at point~$1$) towards a
probability measure~$\nu$, which therefore is a fixed point
of~$\s_\mu$ (see \cite{M1,M2,K1,KP,DL}). In other words, if $\int x\log(x) \,\mu(\dif x)<\log b$
and if $\bigl(W(w)\bigr)_{w\in\alphabet^*}$ is a family of independent
random variables whose probability distribution is $\mu$, then the
non-negative martingale
\begin{equation}\label{Yn}
Y_n=b^{-n}\sum_{w\in\alphabet^n}W(w|_1)W(w|_2)\cdots W(w|_n)
\end{equation}
is uniformly integrable and converges to a random variable~$Y$ whose
probability distribution $\nu$ belongs to $\mathcal{P}_1$ and
satisfies $\s_\mu\nu =\nu$. This means that there exists~$b$ copies
$W(0),\dots ,W(b-1)$ of $W$ and $b$ copies $Y(0),\dots, Y(b-1)$ of~$Y$
such that these~$2b$ random variables are independent and
\begin{equation}\label{eqfonc}
Y=b^{-1}\sum_{k=0}^ {b-1}W(k)Y(k).
\end{equation}
In this case, we denote the measure~$\nu$ by~$\T\mu $. It is natural
to try and iterate~$\T$. But, in general this is not possible because
$\nu=\T\mu $ may not inherit the property $\int x\log (x)\,\nu(\dif
x)<\log b$. So, we have to find a domain stable under the action
of~$\T$. This will be done by imposing conditions on moments.

Indeed, it is easily seen that the sequence $(Y_n)_{n\ge 1}$ defined
in~(\ref{Yn}) remains bounded in~$L^2$ norm if and only if
$\mathbb{E}(W^2)=\moment_2(\mu)<b$, and that in this case
Formula~(\ref{eqfonc}) yields
\begin{equation}\label{m2}
\esp{Y^2}=\frac{b-1}{b-\esp{W^2}}.
\end{equation}
It follows that if $b\ge 3$ and $\esp{W^2}<b-1$, we have $\esp{Y^2}\le
\esp{W^2}$ (the equality holding only if~$W=1$). Therefore, since the
condition $\esp{W^2}<b$ is stronger than $\mathbb{E}(W\log W)<\log b$
when $\esp{W}=1$ (since the function $t\mapsto \log \esp{W^t}$ is
convex), $\T$~is a transformation on the subset of $\mathcal{P}_1$
defined by
$$\mathscr{P}_b= \bigl\{\mu\in\mathcal{P}_1\ :\ 1< \moment_2(\mu)<b-1
\bigr\}.
$$

If $\mu\in {\mathscr P}_b$, due to (\ref{eqfonc}), we can associate
with each $n\ge 0$ a random variable $W_{n+1}$ as well as $2b$
mutually independent random variables $W_n(0),\dots, W_{n}(b-1)$ and
$W_{n+1}(0),\dots,W_{n+1}(b-1)$ such that
\begin{equation}\label{foncn}
W_{n+1}=\frac{1}{b}\sum_{k=0}^{b-1}W_{n}(k)W_{n+1}(k),
\end{equation}
$\T^n\!\mu $ is the probability distribution of $W_n(k)$ for every~$k$
such that $0\le k\le b-1$, and $\T^{n+1}\!\mu $ is the probability
distribution of $W_{n+1}$ and $W_{n+1}(k)$ for every $0\le k\le b-1$.

In Mandelbrot~\cite{M1,M2}, the random variable $Y$ represents
the increment between~0 and~1 of the non-decreasing continuous
function $h$ on $[0,1]$ obtained as the almost sure uniform limit of
the sequence of non-decreasing continuous functions $\phi_n$ defined
by
\begin{equation}\label{mun}
\phi_n(u)=\int_0^u\prod_{k=1}^nW(\tilde t|_k)\, dt,
\end{equation}
where~$\tilde t$ stands for the sequence of digits in the base~$b$
expansion of~$t$ (of course the ambiguity for countably many~$t$'s is
harmless). In other words, for $w\in \alphabet^*$, we have
\begin{equation}\label{increment}
\Delta(\phi,I_w) = b^{-|w|}Y(w) \prod_{1\le j\le |w|} W(w|_j),
\end{equation}
where $Y(w)$ has the same distribution as~$Y$ and is independent of
the variables $W(w|_j)$.

Let us denote by $F(\mu)$ the probability distribution of the
limit~$\phi$, considered as a random continuous function.

We are going to study the dynamical system
$({\mathscr P}_b,\T)$. This will lead to a description of the
asymptotic behavior of $\big (\T^n\!\mu ,F( \T^{ n-1}\!\mu )\big
)_{n\ge 1}$ as $n$ goes to $\infty$.
\medskip

We need some more definitions. For $b\ge 3$, set
$$ w_{2}(b)=\min \left(b-1,
b\frac{b^4-4b^2+12b-8}{b^4+8b^2-12b+4}\right )
$$  
and, for $t$ such that $1<t< w_{2}(b)$,
$$ w_3(b,t)=\frac{b^2}{2}+\frac12\sqrt{\frac{b(b^4-4b^2+12b-8) -
t\bigl(b^4+8b^2-12b+4\bigr )}{b-t}}.
$$
One always has $w_3(b,t)<b^2-1$. 

Also set 
$$\mathscr{D}_b=\Bigl\{\mu\in\mathcal{P}\ :\ \moment_1(\mu)=1,\, 1<
\moment_2(\mu)< w_2(b),\ \text{and\ }\moment_3( \mu)<
w_3\big(b,\moment_2(\mu)\big )\Bigr\}.
$$ 

\begin{theorem}[Central limit theorem]\label{th1}
Suppose $b\ge 3$. Let $\mu\in \mathscr{P}_{b}$, and, for $n\ge 1$,
define\ \ $\displaystyle\sigma_n=\left (\int (x-1)^2 \,\T^{ n}\!\mu
(\dif x)\right )^{1/2}$. Then
\begin{enumerate}
\item
The limit of $(b-1)^{n/2}\sigma_n$ exists and is positive; so
$\lim_{n\to\infty} \T^{n}\!\mu =\delta_1$.
\item If $\mu\in \mathscr{D}_b$, then, $\displaystyle\sup_{n\ge 1}\int
\left (\frac{|x-1|}{\sigma_n}\right )^3 \,\T^{ n}\!\mu (\dif
x)<\infty$.
\item Suppose that there exists $p>2$ such that
$\displaystyle\sup_{n\ge 1}\int \left (\frac{|x-1|}{\sigma_n}\right
)^p \,\T^{ n}\!\mu (\dif x)<\infty.$
Then, if $W_n$ is a variable whose distribution is $\T^{ n}\!\mu$,
$\displaystyle \frac{W_n - 1}{\sigma_n}$ converges in
distribution towards $\mathcal{N}(0,1)$.
\end{enumerate}
\end{theorem}

\begin{theorem}[Functional CLT]\label{th2}
Suppose $b\ge 3$. Let $\mu\in \mathscr{P}_{b}$. Then
\begin{enumerate}
\item The probability distributions $F(\T^{ n}\!\mu )$ weakly
converges towards~$\delta_{\Id}$.
\item Suppose that there exists $p>2$ such that
$$\displaystyle\sup_{n\ge 1}\int \left (\frac{|x-1|}{\sigma_n}\right
)^p \,\T^{ n}\!\mu (\dif x)<\infty.$$
(In particular, this holds if $\mu$ lies in the domain of attraction
$\mathscr{D}_b$.) Then, if $h_n$ is a random function distributed
according to~$F(\T^{n-1}\!\mu )$, the distribution of $\displaystyle
\frac{h_n-\Id}{\sigma_n}$ weakly converges towards the distribution of
the unique continuous Gaussian process $(X_t)_{t\in [0,1]}$, such that
$X(0)=0$ and, for all $j\ge 1$, the covariance matrix $\mat{M}_j$ of
the vector $\bigl(\Delta (X,I_w)\bigr)_{w\in\alphabet^j}$ is given by
$$ \mat{M}_j(w,w')= \begin{cases}
b^{-2j}\bigl(1+(b-1)|w|\bigr)&\text{if } w= w',\\[4pt]
b^{-2j}(b-1)|w\land w'| &\text{otherwise}.
\end{cases}
$$
\end{enumerate}
\end{theorem}

In Section~\ref{sec4}, we will give an alternate construction of this
Gaussian process~$X$: It will be obtained as the almost sure limit of
an additive cascade of normal variables.

\section{Proof of Theorem~\ref{th1}}

Throughout this section and the next one, we assume $b\ge 3$.

\begin{proposition}\label{moments}
If $\mu\in{\mathscr P}_b$ and $\displaystyle \sigma_n^2=\int
(x-1)^2\T^{ n}\!\mu(\dif x)$, the sequence $(b-1)^{n/2}\sigma_n$
converges to $\displaystyle
\sigma_0\sqrt{\frac{b-2}{b-2-\sigma_0^2}}$.
\end{proposition}

\begin{proof}
Equations (\ref{foncn}) and (\ref{m2}) yield
$\displaystyle\esp{W_{n+1}^2}=\frac{b-1}{b-\esp{W_n^2}}$,
from which we get the formula
\begin{equation}\label{ecarttype}
\sigma_{n+1}^2=\esp{W_{n+1}^2}-1=
\frac{\sigma_n^2}{b-1-\sigma_{n}^2},
\end{equation}
which can be written as
\begin{equation*}
\frac{\sigma_{n+1}^2}{b-2-\sigma_{n+1}^2} =
\frac{(b-1)^{-1}\sigma_{n}^2}{b-2-\sigma_{n}^2}.
\end{equation*}
 This yields 
\begin{equation}\label{ecarttype2}
\frac{ \sigma_n^2}{b-2-\sigma_n^2}=\frac{
\sigma_0^2}{b-2-\sigma_0^2}(b-1)^{-n}.
\end{equation}
\end{proof}

\begin{proposition}
If $\mu\in\mathscr{D}_{b}$, then both sequences
$\bigl(\moment_2(\T^{n}\!\mu )\bigr )_{n\ge 1}$ and $\bigl
(\moment_3(\T^{ n}\!\mu )\bigr)_{n\ge 1}$ are non-increasing and
converge to 1 as $n$ goes to $\infty$.
\end{proposition}

\begin{proof}
For $n\ge 0$, we set $u_n=\moment_2(\T^{ n}\!\mu)$ and
$v_n=\moment_3(\T^{ n}\!\mu)$ and deduce from (\ref{foncn}) that, for
all $n\ge 0$, we have
\begin{eqnarray}
u_{n+1} &=& \frac{b-1}{b-u_n} \label{M2}\\
v_{n+1}&=&\frac{(b-1)\bigl
(3u_nu_{n+1}+b-2\bigr)}{b^2-v_n} \label{M3}
\end{eqnarray}
if $u_n<b$ and $v_n<b^2$.

Since $1\le u_0<b-1$, as we already saw it, Equation~(\ref{M2})
implies that $u_n$ decreases, except in the trivial case
$\mu=\delta_1$. Moreover $u_n$ converges towards~1, the stable fixed
point of $t\mapsto (b-1)/(b-t)$.

The conditions $\moment_2(\mu)\le w_2(b)$ and
$\moment_3(\mu)<w_3\big(b,\mathbf{m_2}(\mu)\big )$ are optimal to
ensure that $v_1\le v_0$, and also they impose $v_0<b^2-1$. We
conclude by recursion: if $v_{n+1}\le v_n<b^2$, then we have
\begin{eqnarray*}
v_{n+2}&\le&
\frac{(b-1)\bigl(3u_{n+1}u_{n+2}+b-2\bigr)}{b^2-v_n}\\
&\le& \frac{(b-1)\bigl(3u_nu_{n+1}+b-2\bigr)}{b^2-v_n} =
v_{n+1}.
\end{eqnarray*}
Thus $(v_n)_{n\ge 0}$ is non-increasing and $1\le v_n<b^2-1$, so
we deduce from (\ref{M3}) and the fact that $u_n$ converges to~1 that
$v_n$ converges to the smallest fixed point of the mapping $x\mapsto
(b^2-1)/(b^2-x)$, namely~1.
\end{proof}

\begin{proposition}\label{ordre3}
There exists $C>0$ such that, for $\mu\in \mathscr{D}_b$ and $n\ge 1$,
we have
\begin{equation*}
\big(b^2-\esp{W_n^3}\big)\esp{|Z_{n+1}^3|}\le\\
(b-1)^{3/2} \esp{|Z_{n}^3|}+C
\bigl((\esp{|Z_{n}^3|})^{2/3}+(\esp{|Z_{n}^3|})^{1/3}+1\bigr),
\end{equation*}
where $Z_n=\displaystyle\frac{W_n-1}{\sigma_n}$.
\end{proposition}

\begin{proof}
We use the following simplified notations: $W=W_n$, $Y=W_{n+1}$,
$\sigma_Y$ and~$\sigma_W$ stand for the standard deviations of~$Y$
and~$W$, $Z_Y=\sigma_Y^{-1}\,|Y-1|$, $Z_W=\sigma_W^{-1}\,|W-1|$, and
$r=\sigma_W/{\sigma_Y}$.

Then Equation~(\ref{eqfonc}) becomes\quad $b\,|Y-1|\le
\displaystyle\sum_{i=0}^{b-1}
W(i)\,|Y(i)-1|+\sum_{i=0}^{b-1}|W(i)-1|$, i.e.,
\begin{equation}
b\,Z_Y \le \displaystyle\sum_{i=0}^{b-1} W(i)\,Z_{Y(i)}+
r\sum_{i=0}^{b-1} Z_{W(i)},
\end{equation}
which yields
$$ \big (b^2-\mathbb{E}(W^3)\big )\mathbb{E}(Z^3_Y)\le
r^3\mathbb{E}(Z_W^3) + \sum_{i=0}^3 \binom{3}{i}\,r^i\,T_i,
$$
where
\begin{eqnarray*}
 T_0 &=& 3(b-1) \esp{W^2}\esp{Z_Y^2}\esp{Z_Y}+(b-1)(b-2)
(\esp{Z_Y})^3,\\
T_1 &=& \esp{(W^2Z_W)}\esp{Z_Y^2} + 2(b-1)\esp(WZ_W)(\esp{Z_Y})^2\\
&{}&\quad{} + (b-1)\esp{Z_W}\esp{W^2}\esp{Z_Y^2} +
(b-1)(b-2)\esp{Z_W}(\esp{Z_Y})^2,\\
T_2 &=& \esp{Z_Y}\esp{(WZ_W^2)} + 2(b-1)\esp(WZ_W)\esp{Z_Y}\esp{Z_W}\\
&{}&\quad{} + (b-1)\esp{Z_Y}\esp{Z_W^2} +
(b-1)(b-2)\esp{Z_Y}(\esp{Z_W})^2,\\
T_3 &=& 3(b-1)\esp{Z_W}\esp{Z_W^2} + (b-1)(b-2)(\esp{Z_W})^3.
\end{eqnarray*} 
As, for $X\in\{W,Y\}$ we have $\esp{Z_X}\le (\esp{Z_X^2})^{1/2}=1$,
and $\esp{X^2}<b$, we get the simpler bound
$$ \big (b^2-\esp{W^3}\big )\esp{Z_Y^3}\le
r^3\esp{Z_W^3} + \sum_{i=0}^3 \binom{3}{i}\,r^i\,T'_i,
$$
where
\begin{eqnarray*}
T'_0 &=& (b-1)(4b-5),\\
T'_1 &=& \esp{(W^2Z_W)}+2(b-1)\esp{(WZ_W)}+(b-1)(2b-3),\\
T'_2 &=& \esp{(WZ_W^2)}+2(b-1) \esp{(WZ_W)}+(b-1)^2,\\
T'_3 &=& b^2-1.
\end{eqnarray*}

Since $\esp{W^3}< b^2-1$, the H\"older inequality yields
$$\esp{(W^2Z_W)}\le (\esp{W^3})^{2/3}(\esp{Z_W^3})^{1/3}\le
(b^2-1)^{2/3}(\esp{Z_W^3})^{1/3}$$
and
$$\esp{(WZ_W^2)}\le (\esp{W^3})^{1/3}(\esp{Z_W^3})^{2/3}\le
(b^2-1)^{1/3} (\esp{Z_W^3})^{2/3}.$$
Furthermore, $\esp{(WZ_W)}\le \left( \esp{W^2}
\esp{Z_W^2}\right)^{1/2} \le \sqrt{b-1}$.

We know from~(\ref{ecarttype}) that~$r<\sqrt{b-1}$. Therefore, there
exists a constant~$C>0$ independent of $\mu$ such that
$$ \big(b^2-\esp{W^3}\big)\esp{Z_Y^3} \le (b-1)^{3/2}\esp{Z_W^3} +
C\,\bigl( (\esp{Z_W^3})^{2/3}+(\esp{Z_W^3})^{1/3} + 1\bigr).
$$
\end{proof}

\begin{corollary}
If $\mu\in\mathscr{D}_b$ then\quad $\displaystyle \sup_{n\ge 1}\int
\sigma_n^{-3}\,|x-1|^3 \,\T^{ n}\!\mu (\dif x)<\infty$.
\end{corollary}

\begin{proof} Since $b^2-\esp{W_n^3}$ converges towards $b^2-1$,
and $b^2-1>(b-1)^{3/2}$, the bound in the last proposition yields that
$Z_n$ is bounded in $L^3$.
\end{proof}
\medskip

Recall that we set $\displaystyle Z_n =
\frac{W_n-1}{\sigma_n}$. Equation~(\ref{foncn}) yields
\begin{equation}\label{c6}
Z_{n+1} = \\\frac{1}{b} \sum_{k=0}^{b-1} \left[\sigma_n\,Z_{n}(k)\,
Z_{n+1}(k) + \frac{\sigma_n}{\sigma_{n+1}}\, Z_{n}(k) + Z_{n+1}(k)
\right].
\end{equation} 
 
If we set
$$R_n = \frac{1}{b} \sum_{j=0}^{b-1} Z_{n}(j) Z_{n-1}(j) \sigma_{n-1}
+ \frac{1}{b} \left( \frac{\sigma_{n-1}}{\sigma_n} - \sqrt{b-1}
\right) \sum_{j=0}^{b-1} Z_{n-1}(j),
$$
then Equation~(\ref{c6}) rewrites as
\begin{equation}\label{c7}
Z_{n+1} = R_{n+1} +\frac{\sqrt{b-1}}{b} \sum_{k=0}^{b-1} Z_n(k) +
\frac{1}{b} \sum_{k=0}^{b-1} Z_{n+1}(k).
\end{equation}

We are going to use repeatedly Formula~(\ref{c7}). Fix~$n>1$ and write
\begin{equation}\label{c8}
Z_{n}=Z_n(\epsilon,\epsilon) = R_{n}(\epsilon,\epsilon)
+\frac{\sqrt{b-1}}{b} \sum_{j\in \alphabet} Z_{n-1}(j,0) + \frac{1}{b}
\sum_{j\in \alphabet} Z_{n}(j,1)
\end{equation}
and
\begin{eqnarray*}
Z_{n}(j,1) &=& R_{n}(j,1) +\frac{\sqrt{b-1}}{b} \sum_{k\in \alphabet}
Z_{n-1}(jk,10) + \frac{1}{b} \sum_{k\in \alphabet} Z_{n}(jk,11)\\
Z_{n-1}(j,0) &=& R_{n-1}(j,0) +\frac{\sqrt{b-1}}{b} \sum_{k\in
\alphabet} Z_{n-2}(jk,00) + \frac{1}{b} \sum_{k\in \alphabet}
Z_{n-1}(jk,01).
\end{eqnarray*}
Then Formula~(\ref{c8}) rewrites as
\begin{multline*}
Z_n(\epsilon,\epsilon) = R_n(\epsilon,\epsilon)+ b^{-1}\sum_{j\in
\alphabet} \left( \sqrt{b-1}\,R_{n-1}(j,0) + R_{n}(j,1)\right)+ \ {}\\
b^{-2} \sum_{w\in \alphabet^2}\left( (b-1)\,Z_{n-2}(w,00) +
\sqrt{b-1}\,\bigl(Z_{n-1}(w,01) + Z_{n-1}(w,10)\bigr) +
Z_{n}(w,11)\right),
\end{multline*}
and so on. At last we get $Z_n = T_{1,n} + T_{2,n}$, with
\begin{eqnarray}
\label{T1n} T_{1,n} 
&=& \sum_{k=0}^{n-1} b^{-k} 
\sum_{\substack{m\in\{0,1\}^k\\w\in \alphabet^k}}
(b-1)^{(k-\varsigma(m))/2}R_{n-k+\varsigma(m)}(w,m)\label{r1}\\
T_{2,n} &=& b^{-n}\sum_{\substack{m\in\{0,1\}^n\\w\in\alphabet^n}}
(b-1)^{(n-\varsigma(m))/2} Z_{\varsigma(m)}(w,m),\label{r2}
\end{eqnarray}
where~$\varsigma(m)$ stands for the sum of the components of~$m$.
Moreover, all variables in Equation~(\ref{r2}) are independent, and in
Equation~(\ref{r1}), the variables corresponding to the same~$k$ are
independent.

\begin{proposition}\label{prop1}
We have $\displaystyle \lim_{n\to\infty} \esp{T_{1,n}^2} = 0$, so
$T_{1,n}$ converges in distribution to 0.
\end{proposition}

\begin{proof}  Set $r^2_n = \esp{R_n^2}$. We have
\begin{equation*}
b\, r_n^2 = \sigma_{n-1}^2+\left(
\frac{\sigma_{n-1}}{\sigma_n}-\sqrt{b-1}\right)^2,
\end{equation*}
which together with Formulae~(\ref{ecarttype}) and~(\ref{ecarttype2})
implies that there exists $C>0$ such that $r_n^2\le C(b-1)^{-n}$ for
all $n\ge 1$. By using the independence properties of random variables
in~(\ref{T1n}) as well as the triangle inequality, we obtain
\begin{eqnarray*}
\bigl(\esp{T^2_{1,n}}\bigr)^{1/2}&\leq &\sum_{0\le k< n} b^{-k}
\left( \sum_{0\le j\le k} \binom{k}{j}\,b^k (b-1)^{j}
\,r_{n-j}^2\right)^{1/2}\\
&\le & C\, \sum_{0\le k< n} b^{-k}\left( \sum_{0\le j\le k}
\binom{k}{j}\,b^k (b-1)^{j} \, (b-1)^{j-n}\right)^{1/2}\\
&\le& C\, \sum_{0\le k< n} b^{-k/2} \big ((b-1)^2+1\big )^{k/2}
(b-1)^{-n/2}\\
&=&C\, (b-1)^{-n/2} \sum_{0\le k< n} \left
(\frac{(b-1)^2+1}{b}\right)^{k/2}\\ &=&\mathrm{O}\left( \Big
(1 - \frac{b-2}{b(b-1)}\Big )^{n/2}\right ).
\end{eqnarray*}
\end{proof}

\begin{proposition}\label{prop2}
If there exists $p>2$ such that
$$\sup_{n\ge 1}\int \left (\frac{|x-1|}{\sigma_n}\right )^p \,\T^{
n}\!\mu (\dif x)<\infty,$$
(i.e., $(|Z_n|)_{n\ge 1}$ is bounded in $L^p$), then $T_{2,n}$
converges in distribution to ${\mathcal N} (0,1)$.
\end{proposition}

\begin{proof} 
If $Y$ is a positive random variable, $a$, $p$ and~$\varepsilon$
are positive numbers with~$p>2$, we have
\begin{eqnarray*}
\esp \left(a^2Y^2\,\ind_{\{aY>\varepsilon\}}\right) &\le& a^2 \left(
\esp Y^p\right)^{2/p} \,
\bigl(\mathbb{P}(aY>\varepsilon)\bigr)^{1-2/p}\\
&\le& a^2\left( \esp Y^p\right)^{2/p}\bigl( \varepsilon^{-p} a^{p}
\esp Y^p\bigr)^{1-2/p} = a^p\varepsilon^{2-p} \esp Y^p.
\end{eqnarray*}
So, we have
\begin{multline*}
\sum_{\substack{m\in\{0,1\}^n\\w\in\alphabet^n}}
b^{-2n}(b-1)^{(n-\varsigma(m))} \esp\left(Z_{\varsigma(m)}(w,m)^2
\ind_{\bigl\{b^{-n}(b-1)^{(n-\varsigma(m))/2}
|Z_{\varsigma(m)}(w,m)|> \varepsilon\bigr\}}\right)\\
= \sum_{k=0}^{n} \binom{n}{k}\,b^{n-np}(b-1)^{p(n-k))/2}
\varepsilon^{2-p} \esp |Z_{k}|^p \le \left(
\frac{(b-1)^{p/2}+1}{b^{p-1}}\right)^n\sup_{k\ge 0} \esp |Z_k|^p,
\end{multline*}
and this last expression converges towards~0 as~$n$ goes to~$\infty$.
But, as we have
\begin{equation*}
\esp T_{2,n}^2 = \sum_{\substack{m\in\{0,1\}^n\\w\in\alphabet^n}}
b^{-2n}(b-1)^{n-\varsigma(m)} = \sum_{k=0}^{n}
\binom{n}{k}\,b^{-n}(b-1)^{n-k} = 1.
\end{equation*}
The Lindeberg theorem yields the conclusion.
\end{proof}

\section{Proof of Theorem~\ref{th2}}

We begin by the following observation:
for any real function~$f$ on~$[0,1]$, one has
\begin{equation}\label{holder'}
\omega(f,\delta) \le 2(b-1)\!\sum_{j\ge -\frac{\log \delta}{\log b}}
\sup_{w\in \alphabet^j} \Delta(f,I_w),
\end{equation}
where, $\omega(f,\delta)$ stands for the modulus of continuity of a
function $f$ on $[0,1]$:
\begin{equation*}
\omega(f,\delta)=\sup_{\substack{t,s\in[0,1]\\ |t-s|\le
\delta}}|f(t)-f(s)|.
\end{equation*}

\begin{proposition}\label{tension}
Suppose that $\mu\in{\mathscr P}_b$. If $h_n$ is a random continuous
function distributed according to~$F(\T^{n-1}\!\mu )$, set
$\mathcal{Z}_n=\displaystyle \frac{h_n-\Id}{\sigma_n}$. The
probability distributions of the random continuous functions
$\mathcal{Z}_n$, $n\ge 1$, form a tight sequence.
\end{proposition}

\begin{proof}

By Theorem 7.3 of \cite{Bil}, since $(h_n-\Id)(0)=0$ almost surely for
all $n\ge 1$, it is enough to show that for each positive
$\varepsilon$

\begin{equation}  \label{tightness}
\lim_{\delta\to 0}\limsup_{n\to\infty} \mathbb{P}\big(\omega(
\mathcal{Z}_n,\delta)\ge 2(b-1)\,\varepsilon \big)=0,
\end{equation}

We first establish the following lemma.
\begin{lemma}\label{holder}
Let $\gamma$ and~$H$ be two positive numbers such that
$2H+\gamma-1<0$. Also let $n_0\ge 1$ be such that $\sup_{n\ge
n_0-1}\esp{W_n^2}\le b^\gamma$. For~$j\ge 1$,~$n\ge n_0$ and~$t>0$ we
have
$$ \mathbb{P}\Bigl(\sup_{w\in\alphabet^j}
\Delta\bigl(\mathcal{Z}_n,I_w\bigr)\ge t\,b^{-jH}\Bigr)\le
(b-1)\,t^{-2}(j+1)^3b^{j(2H+\gamma-1)}.$$
\end{lemma}
\begin{proof}
Let $j\ge 1$, $w\in\alphabet^j$ and $n\ge
n_0$. Formula~(\ref{increment}) shows that the increment
$\Delta_n(w)=\Delta (\mathcal{Z}_n,I_w)$ takes the form
\begin{equation}\label{decomp}
\begin{split}
\Delta_n(w) &{}=
b^{-j}\sigma_n^{-1}\left[W_{n}(w)\prod_{k=1}^jW_{n-1}(w|_k)-1\right]\\
&{}= b^{-j}Z_n(w)\prod_{k=1}^jW_{n-1}(w|_k)\\ &\ \qquad+
b^{-j}\sum_{l=1}^j \frac{\sigma_{n-1}}{\sigma_n}Z_{n-1}(w|_l)
\prod_{k=1}^{l-1}W_{n-1}(w|_k).
\end{split}\end{equation}
Consequently,
\begin{multline*}
\mathbb{P}\left(|\Delta_n(w)|\ge t\,b^{-jH} \right)\le \mathbb{P}\left
(b^{-j}Z_n(w)\prod_{k=1}^jW_{n-1}(w|_k)\ge
\frac{t\,b^{-jH}}{j+1}\right )\\+\sum_{l=1}^j\mathbb{P}\left
(b^{-j}\frac{\sigma_{n-1}}{\sigma_n}Z_{n-1}(w|_l)
\prod_{k=1}^{l-1}W_{n-1}(w|_k)\ge \frac{t\,b^{-jH}}{j+1}\right ).\
\end{multline*}
By using the Markov inequality, the equality $\esp{Z_k^2}=1$, and the
fact that $\esp{W_{n-1}^2} \ge 1$, we obtain that each probability in
the previous sum is less than
$$(b-1)\,t^{-2}(j+1)^2 b^{-2(1-H)j} \bigl(\esp{W_{n-1}^2}\bigr)^{j},$$
so that the sum of these probabilities is bounded by
$(b-1)\,t^{-2}(j+1)^3b^{j(\gamma-2(1-H))}$.  Consequently,
\begin{eqnarray*}
\mathbb{P}\left (\exists\ w\in\alphabet^j,\ |\Delta_n(w)|\ge
t\,b^{-jH}\right )&\le& (b-1)\,t^{-2}(j+1)^3b^{j(\gamma-2(1-H)+1)}\\
&=& (b-1)\,t^{-2}(j+1)^3b^{j(2H+\gamma-1)}.
\end{eqnarray*}
\end{proof}

Now, we can continue the proof of Proposition~\ref{tension}.  Fix~$H$,
$\gamma$, and~$n_0$ as in Lemma~\ref{holder}, set $j_\delta = -\log_b
\delta$, and assume that~$n\ge n_0$. Due to~(\ref{holder'}) and
Lemma~\ref{holder}, we have
\begin{eqnarray*}
\bigl\{\omega( \mathcal{Z}_n,\delta)\ge 2(b-1)\,\varepsilon \bigr\}
&\subset& \left\{ \sum_{j\ge j_\delta} \sup_{w\in \alphabet^j}
\Delta(\mathcal{Z}_n,I_w) > \varepsilon\right\}\\
&\subset& \bigcup_{j\ge j_\delta}\left\{ \sup_{w\in \alphabet^j}
\Delta(\mathcal{Z},I_w) >
(1-b^{-H})\,b^{j_0H}\,\varepsilon\,b^{-jH}\right\},
\end{eqnarray*}
so
\begin{equation*}
\mathbb{P}\bigl(\omega( \mathcal{Z}_n,\delta)\ge 2(b-1)\,\varepsilon
\bigr) \le \frac{(b-1)\,b^{-2j_0H}}{(1-b^{-H})^{2}\varepsilon^{2}}
\sum_{j\ge j_0} (j+1)^3 b^{(2H+\gamma-1)\,j}.
\end{equation*}

Consequently,
\begin{equation*}
\lim_{\delta\to 0} \sup_{n\ge n_0} \mathcal{P}\bigl(\omega
(\mathcal{Z}_n,\delta)>2(b-1)\,\varepsilon\bigr) =0.
\end{equation*}
\end{proof}

\begin{proposition}
Suppose that $\mu\in \mathscr{D}_b$. For every $n\ge 1$ let $h_n$ be a
random continuous function whose probability distribution if $F(\T^{
n}\!\mu )$. Fix $j\ge 1$. The probability distribution of the vector
$\left(\Delta (\frac{h_n-\Id}{\sigma_n},I_w)\right)_{w\in\alphabet^j}$
converges, as~$n$ goes to~$\infty$, to that of a Gaussian vector whose
covariance matrix $\mat{M}_j$ is given by
$$\mat{M}_j(w,w')= \begin{cases} b^{-2j}(1+(b-1)|w|)&\text{if}\ w=
w',\\[2pt] b^{-2j}(b-1)|w\land w'| &\text{otherwise}
\end{cases}.
$$
\end{proposition}

\begin{proof}
We use the same notations as in the proof of Lemma~\ref{holder}. Let
$j\ge 1$ and $w\in\alphabet^j$. In the right hand side of
(\ref{decomp}), the random variables $Z_n(w)$ and $Z_{n-1}(w|_l)$,
$1\le l\le j$, are mutually independent and their probability
distribution converge weakly to $\mathcal{N}(0,1)$, while the common
probability distribution of the $W_{n-1}(w|_l)$, $1\le l\le j$,
converges to $\delta_1$, and $\frac{\sigma_{n-1}}{\sigma_n}$ converges
to $\sqrt{b-1}$.

This implies that there exist
$\bigl(\mathcal{N}(v)\bigr)_{v\in\bigcup_{k=1}^j\alphabet^k}$ and
$\bigl(\widetilde{\mathcal{N}}(w)\bigr)_{w\in\alphabet^j}$ two
families of $\mathcal{N}(0,1)$ random variables so that all the random
variables involved in these families are mutually independent, and
\begin{equation}\label{margin}
\lim_{n\to\infty}\bigl(\Delta_n(w)\bigr)_{w\in \alphabet^j}
\stackrel{\scriptscriptstyle\mathrm{dist}}{=}
b^{-j}\left(\widetilde{\mathcal{N}}(w)+ \sqrt{b-1}\sum_{k=1}^j\,
\mathcal{N}(w|_k)\right )_{w\in\alphabet^j}.
\end{equation}
The fact that the vector in the right hand side of (\ref{margin}) is
Gaussian is an immediate consequence of the independence between the
normal laws involved in its definition. The computation of the
covariance matrix is left to the reader.
\end{proof}

\section{The limit process as the limit of an additive
cascade}\label{sec4}

Recall that, if $v\in\alphabet^*$, $[v]$ stands for the cylinder
in~$\alphabet^\omega$ consisting of sequences beginning by~$v$.  Let
$\alphabet^+$ stand for the set of non-empty words on the
alphabet~$\alphabet$.

We are going to show that there exists a finitely additive random
measure~$M$ on $\alphabet^\omega$ satisfying almost surely for all
$w\in \alphabet^+$
\begin{equation}\label{measure}
M([w]) = b^{-j}\left(\zeta(w)+ \sqrt{b-1}\sum_{k=1}^j\,
\xi(w|_k)\right)
\end{equation}
instead of~(\ref{margin}), where the variables
$\bigl(\xi(w)\bigr)_{w\in \alphabet^+}$ are independent with common
distribution~$\mathcal{N}(0,1)$ and the variable $\zeta(w)$ is
$\mathcal{N}(0,1)$ and independent of $\bigl( \xi(w|_j)\bigr)_{1\le
j\le |w|}$.

Indeed, if we set $S(w) = \sum_{k=1}^j\, \xi(w|_k)$, we
should have
\begin{eqnarray*}
b\,\bigl( \zeta(w)+ \sqrt{b-1}\, S(w)\bigr) &=&
\sum_{\ell\in \mathscr{A}} \left(\zeta(w\ell)+
\sqrt{b-1}\sum_{k=1}^{j+1}\, \xi\bigl((w\ell)|_k\bigr)\right)\\
&=& \sum_{\ell\in \mathscr{A}} \left(\zeta(w\ell) +
\sqrt{b-1}\,\xi(w\ell)\right) + b\sqrt{b-1}\,S(w).
\end{eqnarray*}
Iterating this last formula, gives
\begin{equation*}
\zeta(w) = b^{-n} \sum_{v\in \mathscr{A}^n} \zeta(wv) +
\sqrt{b-1}\sum_{j=1}^{n} b^{-j}\sum_{v\in \mathscr{A}^j} \xi(wv).
\end{equation*}
The first term of the right hand side converges to~0 with
probability~1 since its $L^2$ norm is~$b^{-n/2}$. The second term is a
martingale bounded in $L^2$ norm. Therefore its limit, a
$\mathcal{N}(0,1)$ variable, is a.s.~equal to $\zeta (w)$.

Finally, we get a finitely additive Gaussian random measure defined on
the cylinders of $\alphabet^\omega$ by
\begin{equation}\label{measure2}
M([w]) = b^{-|w|} \sqrt{b-1}\,\left(\lim_{n\to\infty}\sum_{v\in
\bigcup_{k=1}^n\alphabet^k} b^{-|v|}\xi(wv) + \sum_{1\le k\le |w|}
\xi(w|_k) \right).
\end{equation}

Then, the limit process of the previous sections can be seen as the
primitive of the projection of~$M$ on $[0,1]$.

Of course~(\ref{measure2}) makes sense even for~$b=2$.

It is easy to compute covariances:
\begin{equation*}
\esp\bigl(M([v])M([w]\bigr) =
\begin{cases} b^{-2|v|}(1+(b-1)|v|)&\text{if}\ w=v,\\[2pt]
(b-1)\, b^{-(|v|+|w|)}|v\land w| &\text{otherwise}.
\end{cases}
\end{equation*}
It is then straighforward to check that, with probability~1, for
all~$\varepsilon>0$ we have $\sup_{v\in\mathscr{A}^n} |M([v])| =
\mathrm{o}(b^{-n(1-\varepsilon)})$. This can be refined, in particular
thanks to the multifractal analysis of the branching random walk
$S(w)= \sum_{1\le j\le |w|}\xi(w|_j)$. In term of the associated
Gaussian process $(X_t)_{t\in [0,1]}$, it is natural to consider for
all $\alpha\in\mathbb{R}$ the sets
$$ \overline{E}_\alpha=\left \{t\in [0,1):
\limsup_{n\to\infty}\frac{\Delta (X,I_n(t))}{n
b^{-n}}=\alpha\sqrt{b-1}\right \},
$$
$$ \underline{E}_\alpha=\left \{t\in [0,1):
\liminf_{n\to\infty}\frac{\Delta (X,I_n(t))}{n
b^{-n}}=\alpha\sqrt{b-1}\right \},
$$
and
$$
E_\alpha=\underline{E}_\alpha\bigcap \overline{E}_\alpha,
$$ where $I_n(t)$ stands for the semi-open to the right $b$-adic
interval of generation $n$ containing $t$.

In the next statement, $\dim\, E$ stands for the Hausdorff dimension
of the set $E$.

\begin{theorem}\label{etudefine}
With probability 1, 
\begin{enumerate}
\item The modulus of continuity of $X$ is a
$\mathrm{O}\bigl(\delta\log(1/\delta)\bigr)$,
\item $X$ does not belong to the Zygmund class,
\item the set $E_0$ contains a set of full Lebesgue measure at each
point of which $X$ is not differentiable,
\item $\displaystyle \dim\, E_\alpha=\dim \,
\underline{E}_\alpha=\dim\, \overline{E}_\alpha
=1-\frac{\alpha^2}{2\log b}$ if $|\alpha|\le\sqrt{2\log b}$, and
$E_\alpha=\emptyset$ if $|\alpha|>\sqrt{2\log b}$. Furthermore,
$E_{-\sqrt{2\log b}}$ and $E_{\sqrt{2\log b}}$ are nonempty.
\end{enumerate}
\end{theorem}

\begin{remark}
We do not know whether there are points in~$E_0$ at which~$X$ is
differentiable. We do not either if the pointwise regularity of~$X$
is~1 everywhere.
\end{remark}

\begin{proof}
For $w\in\alphabet^*$, as above we set
$\zeta(w)=\sqrt{b-1}\lim_{n\to\infty}\sum_{v\in
\bigcup_{k=1}^n\alphabet^k} b^{-|v|}\xi(wv)$.

We have 
$\sum_{n\ge 1}\mathbb{P} \bigl(\exists\ w\in \alphabet^n,\ 
|\zeta(w)|> 2\sqrt{2\log b} \sqrt{n}\bigr)<\infty$, hence,  with
probability~1,
$\sup_{w\in \alphabet^n}|\zeta(w)| = \mathrm{O}(\sqrt{n})$.
\medskip

Also, $\sum_{n\ge 1}\mathbb{P}\bigl( \exists \ w\in \alphabet^n,\
|S(w)|>2n\sqrt{2\log b}\bigr) <\infty$ and, with
probability~1,
$\sup_{w\in \alphabet^n}|S(w)| = \mathrm{O}(n)$.
This yields the property
regarding the modulus of continuity thanks to (\ref{holder'}).
This proves the first assertion.

To see that $X$ is not in the Zygmund class, it is enough to find
$t\in (0,1)$ such that $\displaystyle\limsup_{h\to 0,h\neq
0}\left|\frac{f(t+h)+f(t-h)-2f(h)}{h}\right|=\infty$. Take $t=b^{-1}$
and $h=b^{-n}$. Let $\underline t$ and $\overline t$ stand for the
infinite words $0(b-1)(b-1)\dots(b-1)\cdots$ and
$1(b-1)(b-1)\dots(b-1)\cdots$. We have
\begin{eqnarray*}
\left|\frac{f(t+h)+f(t-h)-2f(h)}{h\sqrt{b-1}}\right| &=& |s(\overline
t_{|n})-s(\underline t_{|n})|\\
&=& \left|\frac{\zeta (\overline t_{|n})-\zeta (\underline
t_{|n})}{\sqrt{b-1}}+S (\overline t_{|n})-S (\underline
t_{|n})\right|.
\end{eqnarray*}
Since $|\zeta (\overline t_{|n})-\zeta (\underline
t_{|n})|=\mathrm{O}(\sqrt{n})$ and since the random walks $S
(\overline t_{|n})$ and $S (\underline t_{|n})$ are independent, the
law of the iterated logarithm yields the desired behavior as
$h=b^{-n}$ goes to~$0$.
\medskip

The fact that $E_0$ contains a set of full Lebesgue measure on which
$X$ is nowhere differentiable is a consequence of the Fubini theorem
combined with the property $|\zeta (\widetilde t_{|n})|=O(\sqrt{n})$
and the law of the iterated logarithm which almost surely holds for
the random walk $(S(\widetilde t_{|n}))_{n\ge 1}$ for each $t\in
[0,1]$.
\medskip

Since, with probability 1, we have $\sup_{w\in \alphabet^n}|\zeta(w)|=
\mathrm{O}(\sqrt{n})$, we only have to take into account the term
$S(w)$ in the asymptotic behavior of
$\displaystyle\frac{\Delta(X,I_w)}{|w| b^{|w|}}$ as $|w|$ tends to
$\infty$.  Thus, in the definition of the sets $\underline{E}_\alpha$,
$\overline{E}_\alpha$ and $E_\alpha$,
$\displaystyle\frac{\Delta(X,I_w)}{\sqrt{b-1} |w| b^{|w|}}$ can be
replaced by $\displaystyle \frac{S(w)}{|w|}$. Then the result is
mainly a consequence of the work~\cite{B1} on the multifractal
analysis of Mandelbrot measures.

To get an upper bound for the Hausdorff dimensions, we set
$$
\beta(q)=\liminf_{n\to\infty}-\frac{1}{n}\log_b\sum_{w\in\alphabet^n}
\exp(qS(w))
$$
for $q\in\mathbb{R}$. Standard large deviation estimates show that
$\displaystyle \dim \, F_\alpha\le \inf_{q\in\mathbb{R}} \frac{-\alpha
q}{\log b }-\beta(q)$ for all $\alpha\in\mathbb{R}$ and
$F\in\{\underline{E},\overline{E},E\}$ (the occurrence of a negative
dimension meaning that the corresponding set is empty). Also, using
the fact that $\beta(q)$ is the supremum of those numbers $t$ such
that $\limsup_{n\to\infty} b^{nt}\sum_{w\in\alphabet^n}
\exp(qS(w))<\infty$ yields
$$\beta(q)\ge\lim_{n\to \infty} \displaystyle-\frac{1}{n}\log_b
\mathbb{E}\sum_{w\in\alphabet^n} \exp(qS(w))= -1-\frac{q^2}{2\log b }.
$$
Since both sides of this inequality are concave functions, we actually
have, with probability 1, $\beta(q)\ge -1-q^2/2\log b $ for all
$q\in\mathbb{R}$. Consequently, the upper bound for the dimension used
with $\alpha=-\beta'(q)\log b =q$ yields, with probability 1, $\dim \,
F_q\le 1-q^2/2$ for all $q\in [-\sqrt{2\log b },\sqrt{2\log b }]$ and
$F_q=\emptyset$ if $|q|>\sqrt{2\log b }$.
\smallskip

For the lower bounds, we only have to consider the sets $E_\alpha$.

If $q\in [-\sqrt{2\log b },\sqrt{2\log b }]$, let $\phi_q$ be the
non-decreasing continuous function associated with the family $\big
(W_q(w)=\exp(qW(w)-q^2/2)\big )_{w\in\alphabet^+}$ as $\phi$ was with
$\big (W(w)\big )_{w\in\alphabet^+}$ in Section~\ref{dynam}. We learn
from \cite{B1} that, with probability 1, all the functions $\phi_q$,
$q\in (-\sqrt{2\log b },\sqrt{2\log b })$ are simultaneously defined;
their derivatives (in the sense of distributions) are positive
measures denoted by $\mu_q$.  Then, computations very similar to those
used to perform the multifractal analysis of $\mu_1$ in \cite{B1} show
that, with probability 1, for all $q\in (-\sqrt{2\log b },\sqrt{2\log
b })$ the dimension of $\mu_q$ is $1-q^2/2\log b $ and $\mu_q(E_q)>0$.

For $q\in\{-\sqrt{2\log b },\sqrt{2\log b }\}$ it turns out (see
\cite{Liu,B1}) that the formula
$$ \mu_q(I_w)=\lim_{p\to\infty} -\sum_{v\in
\alphabet^p}\Pi(wv)\log\Pi(wv),
$$
where
$$\Pi(wv)=b^{-|w|+p}\prod_{k=1}^{|w|}W_q(w_{|k})
\prod_{j=1}^pW_q(wv_{|j}),$$
defines almost surely a positive measure carried by $E_q$.
\end{proof}

\section{Other random Gaussian measures and processes}

This time~$b\ge2$, $\bigl( \xi(w)\bigr)_{w\in \mathscr{A}^+}$ is a
sequence of independent $\mathcal{N}(0,1)$ variables, and
$\bigl(\alpha(w)\bigr)_{w\in \mathscr{A}^+}$ and
$\bigl(\beta(w)\bigr)_{w\in \mathscr{A}^+}$ are sequences of numbers
subject to the conditions
\begin{eqnarray*}
&{}& \alpha(w) = \sum_{\ell\in \mathscr{A}}
  \alpha(w\ell),\\
&{}& \sum_{v\in \mathscr{A}^+} |\alpha(wv)|^p |\beta(wv)|^p <
\infty\quad \mbox{for some $p\in (1,2]$}.
\end{eqnarray*}

Then, for all $w\in \mathscr{A}^+$, the martingale
$\sum_{v\in\mathscr{A}^n} \alpha(wv)\,\beta(wv)\,\xi(wv)$ is bounded
$L^p$ norm (if~$p<2$ this uses an inequality from~\cite{BE}) and the
formula
\begin{equation}\label{generalmeasure}
M([w]) = \lim_{n\to\infty}\sum_{v\in \bigcup_{k=1}^n\alphabet^k}
\alpha(wv)\,\beta(wv)\,\xi(wv) + \alpha(w) \sum_{1\le j\le |w|}
\beta(w|_j)\,\xi(w|_j)
\end{equation}
almost surely defines a random measure which generalizes the one
considered in the previous sections. Here again, the primitive of the
projection on $[0,1]$ of this measure defines a continuous process,
which is Gaussian if $p=2$.

\begin{remark}
The hypotheses under which this last construction can be performed can
be relaxed: if the random variables $\xi(w)$, $w\in \mathscr{A}^+$,
are independent, centered, and $\sum_{v\in \mathscr{A}^+}
|\alpha(wv)|^p |\beta(wv)|^p \mathbb{E}(|\xi(wv)|^p)< \infty$, Formula
(\ref{generalmeasure}) still yields a random measure.

The fine study of the associate process as well as some improvement of
Theorem~\ref{etudefine} will be achieved in a further work.
\end{remark}

\end{document}